\documentclass[letterpaper,11pt]{article}
\usepackage{caption}
\usepackage{graphicx}
\RequirePackage{amsthm,amsmath,amssymb}
\RequirePackage[numbers]{natbib}
\usepackage[ruled,vlined,linesnumbered]{algorithm2e}
\usepackage[usenames,dvipsnames]{color}
\usepackage{ifpdf}
\ifpdf    
\RequirePackage[colorlinks,citecolor=blue,urlcolor=blue]{hyperref}
\else    
\RequirePackage[colorlinks,citecolor=blue,urlcolor=blue,hypertex]{hyperref}
\fi

\usepackage{courier}
\usepackage[export]{adjustbox}

\usepackage[framemethod=tikz]{mdframed}

\newcommand{\argmax}[1]{\underset{#1}{\operatorname{argmax}}\;}

\RequirePackage{hypernat}
\RequirePackage{graphicx,color}
\usepackage{algorithmic}

\newenvironment{Proof}{\begin{proof}}{\end{proof}}
\newtheorem{theorem}{Theorem}

\newtheorem{cor}[theorem]{Corollary}

\newtheorem{lemma}[theorem]{Lemma}

\newtheorem{proposition}[theorem]{Proposition}

\numberwithin{equation}{section}
\numberwithin{theorem}{section}

\newcommand{\comment}[1]{}

\newcommand{\cN}{{\cal{N}}}
\newcommand{\sbeta}{{\sqrt{\beta}}}
\newcommand{\R}{{\mathbb{R}}}

\newcommand{\npk}{(n,p,k)}

\newcommand{\eps}{{\varepsilon}}
\newcommand{\SigHat}{{\Sigma}}
\newcommand{\x}{{\bf x}}

\newcommand{\iid}{i.i.d.}

\newcommand{\rr}{{\bf r}}
\newcommand{\es}{{\bf s}}
\newcommand{\uu}{{\bf u}}
\newcommand{\vv}{{\bf v}}

\newcommand{\xii}{\boldsymbol\xi}

\newcommand\cS{\mathcal{S}}
\newcommand\cI{\mathcal{I}}
\newcommand\cJ{\mathcal{J}}
\usepackage{mathtools}
\newcommand{\onep}{{\{1,\ldots,p\}}}

\DeclareMathOperator{\pr}{Pr}

\title{\Large A greedy anytime algorithm for sparse PCA
\author{Guy Holtzman
        \thanks{Ben-Gurion Unversity of the Negev, Beer-Sheva, Israel}\\
        \and Adam Soffer
        \thanks{Ben-Gurion Unversity of the Negev, Beer-Sheva, Israel}\\
        \and Dan Vilenchik
        \thanks{Ben-Gurion Unversity of the Negev, Beer-Sheva, Israel.}}
}

\begin{document}
\maketitle

\begin{abstract}%

The taxing computational effort that is involved in solving some high-dimensional statistical problems, in particular problems involving non-convex optimization, has popularized the development and analysis of  algorithms that run efficiently (polynomial-time) but with no general guarantee on statistical consistency.  In light of the ever-increasing compute power and decreasing costs, a more useful characterization of algorithms is by their ability to  calibrate the invested computational effort with various characteristics of the input at hand and with the available computational resources. For example, design an algorithm that always guarantees statistical consistency of its output by increasing the running time as the SNR weakens.

\vspace{0.1cm}

We propose a new greedy algorithm for the $\ell_0$-sparse PCA problem which supports the calibration principle. We provide both a rigorous analysis of our algorithm in the spiked covariance model, as well as simulation results and comparison with other existing methods.  Our findings show that our algorithm recovers the spike in SNR regimes where all polynomial-time algorithms fail while running in a reasonable parallel-time on a cluster.
\end{abstract}


\section{Introduction}

 Principal components analysis (PCA) is the mainstay of modern machine learning and statistical inference, with a wide range of applications involving multivariate data, in both science and engineering~\citep{Anderson84,JoliffePCA}. 
The application of PCA to high-dimensional data, where features are plentiful (large $p$) but samples are relatively scarce (small $n$) suffers from two major limitations:
(1) the principal components are typically a linear combination of
all features, which hinders their interpretation and subsequent use;
and (2) while PCA is consistent in the classical setting ($p$ is fixed and $n$ goes to infinity) \cite{Anderson84,Muirhead1982}, it is generally inconsistent in high-dimensions ~\cite{Johnstone01,Paul07,BickelLevina06,Nadler08,Johnstone.Lu2009Consistency}. This phenomenon is more generally known as the ``Curse of Dimensionality"~\cite{Donoho00}.

The lack of consistency and interpretability in the high-dimensional setting encouraged researchers to design regularized estimation schemes, where additional structural information on the parameters describing the statistical models is assumed and in particular various sparsity assumptions. 
One such popular model is the $\ell_0$-sparse PCA, or $k$-sparse PCA as we call it from now on. The input to $k$-sparse PCA is a pair $(X,k)$, where $X$ is an $n \times p$ design matrix and $k$ the desired sparsity level. The goal is to find a unit vector $\vv$ that has at most $k$ non-zero entries, a $k$-sparse vector, such that the variance of $X$ in $\vv$'s direction is maximal.

While  standard (non-restricted) PCA can be efficiently solved by computing the eigenvectors of a symmetric matrix, its $k$-sparse variant has a discrete combinatorial flavor that turns it NP-hard~\cite{Natarajan1995}. Nevertheless, computationally efficient heuristics for $k$-sparse PCA were proposed and analyzed under various assumptions on the distribution of $X$ and the parameters $n,p$ and $k$, e.g.~\cite{Johnstone.Lu2009Consistency,AminiWain09,AspremontEtAlSDP07,krauthgamer2015,Deshpande:2016}.

The performance of all the aforecited algorithms features  a rather undesirable phase-transition behavior (at least on the benchmark distribution that was studied in each paper). Each algorithm $A$ performs well up to a certain SNR threshold $\tau_A$,  and its performance deteriorates as the SNR drops below $\tau_A$ (we formalize the notion of  signal-to-noise-ratio in Section \ref{sec:spikedmodel}). Such a threshold behavior is expected in a worse-case setting, as the problem is NP-hard. However, the results of  \cite{RigolletClique,RigolletCliqueCOLT,krauthgamer2015,Bresler} suggest that the threshold behavior might persist even in the average-case setting, as long as the algorithms belong to the polynomial-time family.

\bigskip

Throughout we let $\vv^* \in \R^p$ denote the solution of the $k$-sparse PCA problem and $\cI^* \subseteq \onep$ the support of $\vv^*$.  We denote by $\SigHat = \tfrac{1}{n} X^TX$ the sample covariance matrix, assuming $X$ is centered. In what follows, we consider the equivalent problem of finding the support set  $\cI^*$ rather than $\vv^*$.

\section{Our contribution} \label{sec:contrib}
Anytime algorithms provide the ability to achieve results of better quality in return for running time~\cite{Zilberstein_1996}. We implement this philosophy in the context of the sparse PCA problem. We propose a new algorithm that consists of a tunable parameter that allows increasing the running time as the SNR weakens. Thus the algorithm maintains a steady success rate and avoids the aforementioned threshold behavior. If necessary, the algorithm invests super polynomial-time.

To understand our algorithm, it will be convenient to reformulate $k$-sparse PCA as follows. For a fixed symmetric matrix $\SigHat \in \R^{p \times p}$, define  $f_{\lambda_1}^{(\SigHat)}:2^{\onep} \to \R$ by  $f_{\lambda_1}^{(\SigHat)}(\cS)=\lambda_{1}(\Sigma_\cS),$ where $\SigHat_\cS$ is the principal submatrix of $\SigHat$ corresponding to the variables in $\cS$ and $\lambda_1$ is its largest eigenvalue. For simplicity we  abbreviate $f_{\lambda_1}^{(\SigHat)}$ by $f_{\lambda_1}$ as $\SigHat$ is always clear from the context. The $k$-sparse PCA problem is the solution of
\begin{equation}\label{eq:defOfSPCA}
  \cI^* = \argmax{\substack{\cS\subseteq \onep \\ |\cS|=k}} f_{\lambda_1}(\cS).
\end{equation}

Our algorithm is composed of two routines. The first, which we call GreedySprasePCA, receives a real valued function $f:2^{\onep} \to \R$ (for example $f_{\lambda_1}$), a seed $\cS^* \subseteq \onep$ of size $k^* \le k$, and a solution size $k$. It greedily completes $\cS^*$ to a candidate solution of $k$-sparse PCA.

\bigskip
\begin{mdframed}
\noindent{\bf{GreedySPCA}($f,\cS^*,k):$}
\begin{algorithmic}[1]
    \STATE $k^* \leftarrow |\cS^*|$
    \FORALL{$i \in \onep \setminus \cS^*$}
    \STATE $a_i \leftarrow f(\cS^* \cup \{i\})$
    \ENDFOR
    \STATE sort the $a_i$'s as $a_{i_1} \ge a_{i_2} \ge ... \ge a_{i_{p-k^*}}$
    \RETURN $\cS^* \cup \{ i_1,\ldots,i_{k-k^*}\}$
\end{algorithmic}
\end{mdframed}

The next routine, SeedSparsePCA (SSPCA for short) enumerates over all possible seeds of a given size $k^*$, completes each one using GreedySPCA, and returns the ``best" solution.

\begin{mdframed}
\noindent{\bf{SSPCA}$(f_1,f_2,k,k^*):$}
\begin{algorithmic}[1]
    \FORALL{seeds $\cS^* \subseteq \{1,\ldots,p\}$ of size $k^*$}
    \STATE $\cS^{(\cS^*)} \leftarrow$ GreedySPCA$(f_1,\cS^*,k)$
    \ENDFOR
    \RETURN $\argmax{\cS^*}f_2(\cS^{(\cS^*)})$
\end{algorithmic}
\end{mdframed}
Note that SSPCA is actually a family of algorithms, depending on the hyper-parameters $f_1,f_2$. We shortly discuss the choice of these parameters.

The running time of SSPCA is $\binom{p}{k^*}O(pk^*+k \log k)$. By varying $k^*$ one obtains a hierarchy of algorithms, which for our choice of $f_1,f_2$ ranges from Diagonal Thresholding \cite{Johnstone.Lu2009Consistency} (for $k^*=0$) up to the naive exhaustive search (for $k^*=k$). The hierarchy realizes the anytime principle.

Another attractive feature of SSPCA, especially for practitioners, is the fact that it is completely white-box with only one simple tunable parameter, $k^*$. Furthermore, SSPCA can  easily be parallelized and run in a multi-core cluster environment. The code  we share is  written in that way.

\vspace{0.1cm}

The following two conditions are sufficient for SSPCA$(f_1,f_2,k,k^*)$ to recover at least $(\delta-\xi)$-fraction of $\cI^*$, for two numbers $\delta,\xi \in [0,1].$
\begin{itemize}
    \item[C1.] There exists a {\em golden seed} $\cS^*$ of size $k^*$ such that GreedySPCA$(f_1,\cS^*,k)$ outputs a set $\cI$ satisfying  $|\cI \cap \cI^*| \ge \delta k$.
    \item[C2.] $\SigHat$ is {\em $\xi$-separable} with respect to $f_2$. Namely, for every two sets $\cI,\cJ$ of size $k$, if $|\cI \cap \cI^*| - |\cJ \cap \cI^*| \ge \xi k$ then $f_2(\cI) > f_2({\cJ})$.
\end{itemize}
For C1 and C2 to be meaningful, one should think of golden seeds with $\delta$ close to 1, and $f_2$-separability with $\xi$ close to 0. Proposition \ref{clm:cond} formally asserts the sufficiency of these conditions.

\vspace{0.2cm}

\begin{figure}[!t]
\centering

\adjustimage{width=17cm,center}{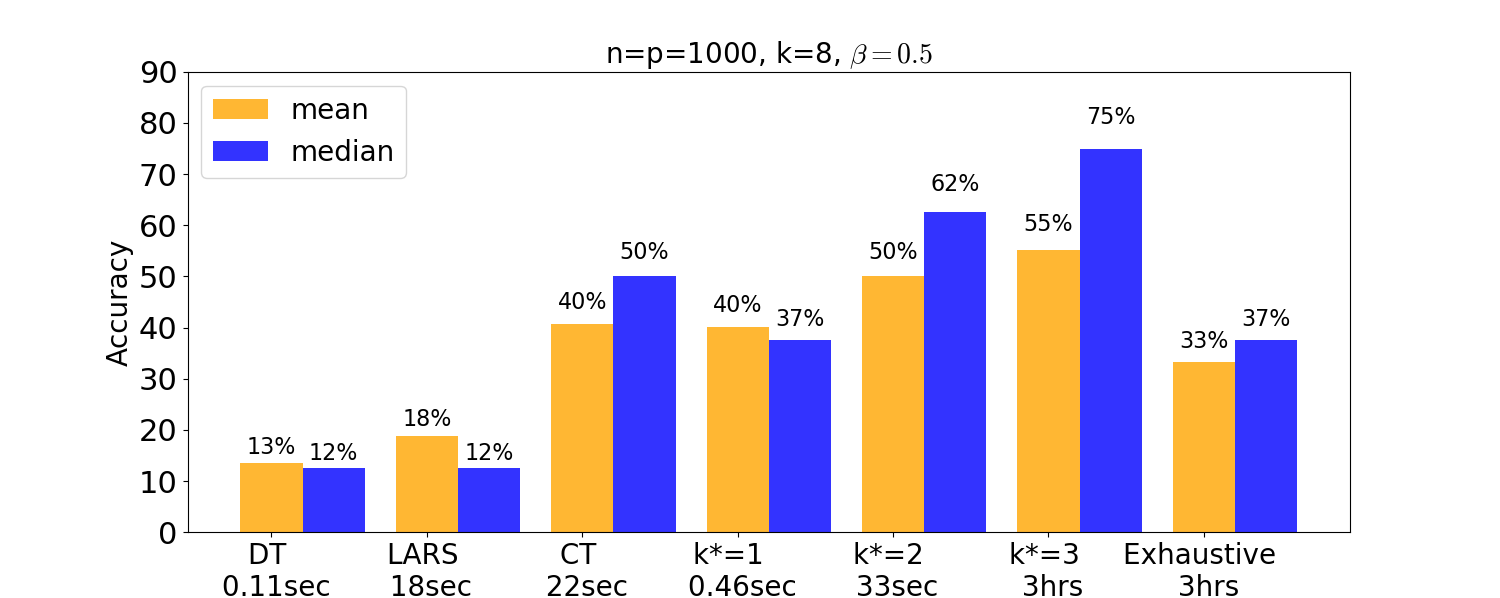}
\caption{\small The plot portrays the accuracy averaged over 25 executions in the uniform unbiased spiked covariance model (USPCA), parametrized to suit a weak SNR regime ($n=p=1000, k=8,\beta=0.5$). The compared algorithms are SSPCA with various seed sizes,  Diagonal Thresholding (DT) \cite{Johnstone.Lu2009Consistency}, Covariance Thresholding (CT) \cite{BickelLevina06}, LARS regression \cite{Zou04sparseprincipal} and a naive exhaustive search that was allowed the same running time as SSPCA with $k^*=3$. Full details of the executions are given in Section \ref{sec:simulation}. The average running time is stated below each algorithm name. The reported running time of SSPCA and the preempted exhaustive search is a parallel-time using 90 {\texttt{Intel Xeon Processor E7-4850 v4 (40M Cache, 2.10 GHz)}} cores.}
\label{fig:comparison}
\end{figure}

The definition of $k$-sparse PCA in Eq.~\eqref{eq:defOfSPCA} suggests the choice $f_2=f_{\lambda_1}$, which is indeed what we chose. For the rigorous analysis, we chose $f_1 = f_{avg}(\cS) =\tfrac{1}{|\cS|}\sum_{i,j \in \cS} \SigHat_{i,j}$, namely the average row sum in $\SigHat_\cS$.  Note that for every $\cS$, $f_{\lambda_1}(\cS) \ge f_{avg}(\cS)$ by plugging the characteristic vector of $\cS$ in the Rayleigh-quotient definition of $\lambda_1$. In the simulation part, we experimented with other functions as well. Details in Section \ref{sec:simulation}.

We analyze SSPCA rigorously in the well-known spiked covariance model, which is formally defined in Section \ref{sec:spikedmodel}.
Theorem \ref{thm:main} establishes the scaling of $k^*$ as a function of the parameters $(n,p,k)$, for which condition C1 holds with $\delta=1$. Theorem \ref{thm:spectral} and Corollary \ref{cor:xi} explicate the gap parameter $\xi$ in condition C2 as a function of $(n,p,k)$, from which the regime where $\xi=o(1)$ is obtained. Together they guarantee the recovery of $(1-o(1))$-fraction of $\cI^*$, up to the information limit. Our results are asymptotic, namely, they hold with probability (w.p.) tending to 1 as the parameters of the problem ($n,p,k$) go to infinity. The probability is taken only over the choice of the design matrix $X$.

\vspace{0.2cm}

Figure \ref{fig:comparison} summarizes simulations that show how our approach implements the anytime paradigm: increasing $k^*$ (and subsequently the running time of SSPCA) translates to the desired increase in the solution quality. We further compared the performance of SSPCA when allowed ``polynomial-time" ($k^*=1,2$) to three popular polynomial-time algorithms. Figure \ref{fig:comparison} shows that SSPCA is better than all three. Finally, we show that SSPCA with $k^*=3$ outperforms the naive exhaustive search when both are running for the same amount of time. This suggests the following rule-of-thumb: Rather than trying to nail down $\cI^*$ by going over as many sets of size $k$ as possible using time budget $T$, we suggest to skim through a larger number of smaller sets of size $k^*<k$, completing each one in a greedy manner (GreedySPCA).

Finally, let us  mention that independently of this result,  a different anytime algorithm for $k$-sparse PCA was obtained in~\cite{ding2019subexponentialtime}. The algorithm also employs a controlled exhaustive search part, but the overall algorithmic approach is different. The algorithm was rigorously analyzed in the spiked covariance model as well and, interestingly, both algorithms have the same asymptotic run-time, namely the same requirement on the seed size. In fact, \cite{ding2019subexponentialtime} provide evidence that this requirement is tight.


\section{The Spiked Covariance Model}\label{sec:spikedmodel}
The spiked covariance model was suggested by Johnstone in 2001~\cite{Johnstone01}  to model a combined effect of a low-dimensional signal buried in high-dimensional noise. In this paper, we consider the  Gaussian case with a single spike, where the population covariance matrix is $\Sigma = \beta\vv^*{\vv^*}^T + I_p$.
The parameter $\beta \ge 0 $ is the signal strength,
\(\vv^* \in \R^p\) is the planted spike assumed to be a $k$-sparse unit-length vector. The algorithmic task is to recover $\cI^*$, the support of $\vv^*$, given $n$ iid samples $\x_1,\ldots,\x_n$ from  $\cN(0,\Sigma)$. The SNR is governed both by $
\beta$ and $k$. The larger $\beta$ and the smaller $k$ the stronger the SNR and the easier the task.

There is an extensive literature on efficient algorithms for sparse PCA, which were rigorously analyzed under various variants of the model just described, e.g.~\cite{AminiWain09, CaiMaWu13, SHEN2013317, Johnstone.Lu2009Consistency, Deshpande:2016,krauthgamer2015,wang2016}. All these algorithms succeed in the regime where the sparsity level satisfies $k =\tilde{O}(\sqrt{\beta^2 n})$ (we use the $\tilde O(\cdot)$ notation in the common way, namely logarithmic factors are ignored). The best sparsity asymptotically is achieved for example by Covariance Thresholding (CT) \cite{Deshpande:2016}, remaining consistent up to $k_0 \asymp \sqrt{\beta^2 n}$ (the notation $f \asymp g$ stands for $f/g  \to c$ for some constant $c > 0$).

It was further shown that under the planted clique hardness assumption there is no polynomial-time algorithm that asymptotically beats $k_0$ \cite{Bresler,RigolletClique,RigolletCliqueCOLT}. Even without the planted clique assumption, \cite{krauthgamer2015} show that the SDP relaxation suggested by \cite{AspremontEtAlSDP07} and analyzed by \cite{AminiWain09}, fails to recover $\vv^*$ beyond $k_0$. The threshold $k_0$ is commonly referred to as the computational threshold, which we denote from now on by  $k_{comp}$. We informally call the regime $k >> k_{comp}$ the weak SNR regime, and $k\le k_{comp}$ the strong SNR regime.
Finally, an information-limit was proven for $k \ge k_{info} \asymp \beta^2 n / \log p$~\cite{AminiWain09,RigolletClique,CMW15,Samworth14}, and a matching algorithmic result for the naive exhaustive search ~\cite{minmax-pmlr-v22-vu12,CaiMaWu13,RigolletCliqueCOLT,pmlr-v75-brennan18a}.

\vspace{0.1cm}

While the boundaries between the different SNR regimes are well understood,
at least asymptotically, the following question remains open:

\vspace{0.1cm}

{\emph{ Question: what is the (minimal) computational complexity required to find the support of $\vv^*$ in the weak SNR regime, namely when $k_{comp} \le k \le k_{info}$?}}

\vspace{0.1cm}

\noindent The analysis of SSPCA  in the next sections provides an answer.

\section{Results}
Our results refer to the following choice of hyper-parameters for SSPCA:  $$f_1=f_{avg}(\cS) =\tfrac{1}{|\cS|}\sum_{i,j \in \cS} \SigHat_{i,j},\qquad  f_2=f_{\lambda_1}(\cS)=\lambda_1(\SigHat_\cS).$$ Furthermore,  we assume the {\em uniform biased sparse PCA model} (UBSPCA), namely non-zero entries of $\vv^*$ are all equal to $1/\sqrt{k}$. In the simulation part we lift the same-sign restriction and use the {\em uniform unbiased sparse PCA model} (USPCA), where entries equal $\pm 1/\sqrt{k}$.

\vspace{0.1cm}

The next proposition asserts the sufficiency of conditions C1, C2.

\begin{proposition}[$\delta,\xi-$Sufficient Conditions] \label{clm:cond}
If $\SigHat$ is $\xi$-separable (Condition C2) and there exists a golden seed $\cS_0$ of size $k^*$ such that GreedySPCA$(X,k,\cS_0)$ outputs a set $\cI_0$ satisfying  $|\cI_0 \cap \cI^*| \ge \delta k$ (Condition C1) then SSPCA outputs a set $\cI$ satisfying $|\cI \cap \cI^*| \ge (\delta-\xi)k$.
\end{proposition}
The proof of Proposition \ref{clm:cond} is given in Section \ref{sec:proof_cond}. Our next theorem  establishes the scaling of $k^*$ for the existence of a seed from which $\cI^*$ is recovered exactly (condition C1 with $\delta=1$).

\begin{theorem}[Golden seed]\label{thm:main} Let $\SigHat$ be distributed according to the UBSPCA model. Assume that $p/n\to
c \ge 0$, and $k\le \frac{n\cdot \min\{\beta^2,\beta\}}{C\log n}$ for a sufficiently large universal constant $C$. If
\begin{equation}\label{eq:kstar}
    k^* \ge \Bigl\lfloor \frac{C k^2 \log n}{\beta^2 n} \Bigr\rfloor
\end{equation}
then w.p. tending to 1 as $\npk \to \infty$ there exists a seed $\cS^* \subseteq \cI^*$ of size at most $k^*$ for which
the output of GreedySPCA$(f_{avg},\cS^*,k)$ is $\cI^*$.
\end{theorem}
Theorem \ref{thm:main} is proven in Section \ref{sec:proofMain}. The next theorem provides a general spectral separability property of spiked covariance matrices. It implies condition C2 as an immediate corollary.
\begin{theorem}[Spectral Separation] \label{thm:spectral}
Let $\SigHat$ be distributed according to the UBSPCA model with $p/n\to c > 0$. Set $\Gamma=C\left(\frac{(1+\beta) k \log n}{n}\right)^{0.5}$ for a suitably chosen constant $C$. With probability tending to 1 as $\npk \to \infty$, for every $\delta \in [0,1]$ and for every set $\cI \subseteq \onep$ of size $k$ that satisfies $| \cI \cap \cI^*|=\delta k$,

$$
\lambda_1(\SigHat_\cI) \in \left[1+\delta\beta -\Gamma, 1+\delta\beta + \Gamma+ \frac{\beta}{k}\right].
$$
\end{theorem}

\begin{cor}[Condition C2]\label{cor:xi} Under the conditions of Theorem \ref{thm:spectral}, for every two sets $\cI,\cJ$ of size $k$, if $|\cI \cap \cI^*| - |\cJ \cap \cI^*| > \xi k$ then $\lambda_1(\SigHat_{\cI}) > \lambda_1(\SigHat_{\cJ})$, for $\xi$ that satisfies
\begin{equation}\label{eq:xi}
  \xi = \tfrac{1}{k}+O\left(\sqrt{\frac{(1+\beta)k}{k_{info}}}\right),
\end{equation}
where $k_{info} \asymp \beta ^2 n /\log p$ was defined in Section. \ref{sec:spikedmodel}
\end{cor}
Theorem \ref{thm:spectral} is proven in Section \ref{sec:ProofSpectral} and Corollary \ref{cor:xi} is proven in Section \ref{sec:proofOfCor}.
Corollary \ref{cor:xi} with $\cI = \cI^*$ was already proven for example in \cite{minmax-pmlr-v22-vu12,CaiMaWu13,RigolletCliqueCOLT,pmlr-v75-brennan18a} and in a more general sparse PCA model. Corollary \ref{cor:xi} adds that even if the seed suffices to recover only part of $\cI^*$, which might as well be the case in practice (finite problem size), SSPCA will nevertheless pick up this information. This point is demonstrated in Figure \ref{fig:comparison}, where all executions of SSPCA end up with partial recovery.

\vspace{0.1cm}

The next theorem concludes our result and is an immediate corollary of all the statements in this section.
\begin{theorem} \label{thm:conclusion}
Under the  conditions of Theorem \ref{thm:main}, if $k^*$ satisfies the lower bound in Eq.~\eqref{eq:kstar}, then w.p. tending to 1 as $(n,p,k) \to \infty$, SSPCA$(f_{avg},f_{\lambda_1},k,k^*)$ recovers at least $(1-O(\sqrt{\alpha})-1/k)$-fraction of $\cI^*$, where $\alpha=(1+\beta)k/k_{info}$.
\end{theorem}

For example, in the regime where DT and SDP fully recover $\cI^*$, i.e. $k=O(k_{comp}/\sqrt{\log p})$, SSPCA requires a seed of
size 0 to recover $\cI^*$ and runs in time $O(p\log p)$. When $k \asymp k_{comp}$ the seed size is $k^*=O(\log n)$ and the running time is quasi-polynomial, $p^{O(\log n)}$. Simulations suggest that up to the computational threshold,  even for a fairly large problem size ($n=p=20,000$), it suffices to choose $k^*=1$   to recover $\cI^*$ exactly (see Figure \ref{fig:boundaries}). In the weak SNR regime, i.e. $k=n^{0.5+\eps}$, the seed size scales as $n^{2\eps}\log n$. This means that the computational effort is  $\exp\{n^{2\eps}\log n\}$, which is  exponential in $(k/\sqrt{n})^2$ (square the excess above the computational threshold)  rather than in $k$ itself (the naive exhaustive search approach). The results in \cite{ding2019subexponentialtime}[Thm~2.14] provide rigorous evidence that the exact scaling that we obtained for $k^*$ in Eq.~\eqref{eq:kstar} is asymptotically optimal.

\vspace{0.1cm}

Simulation (Section \ref{sec:simulation}) suggests that SSPCA succeeds also when the UBSPCA assumption is  relaxed, namely the same-sign assumption is lifted. In this case, the best performance is achieved  when the hyper-parameter $f_{avg}$ is replaced with $f_{\ell_1}$, which is the average row $\ell_1$-norm rather than the average row sum. The parameter $f_2=f_{\lambda_1}$ remains the same.


 \section{Proof of Proposition \ref{clm:cond}}\label{sec:proof_cond}
Suppose by contradiction that conditions C1 and C2 hold but SSPCA outputs a set $\cJ$ for which $|\cJ \cap \cI^*| < (\delta-\xi)k$. Consider a point in the execution of SSPCA where a golden seed $\cS_0$ is explored. By Condition C1, GreedySPCA completes $\cS_0$ to a set $\cI$ satisfying $|\cI \cap \cI^*| \ge \delta k$. The latter together with the contradiction assumption give  $|\cI \cap \cI^*| - |\cJ \cap \cI^*| > \delta k - (\delta-\xi)k  = \xi k$. In this case C2 guarantees that $\lambda_1(\SigHat_{\cI}) > \lambda_1(\SigHat_{\cJ})$. Therefore the last line of SSPCA ensures that $\cJ$ cannot be the output of the algorithm.


\section{Proof of Theorem \ref{thm:main}}\label{sec:proofMain}
For convenience, let us assumes w.l.o.g. that the support of $\vv^*$ is the first $k$ variables, namely $\cI^*=\{1,\ldots,k\}$. Our candidate for a golden seed is any subset $\cS^*$ of $\cI^*$. For concreteness we fix $\cS^* =  \{1,\ldots,k^*\}$.  We show that when GreedySPCA is called with this subset, then the $k-k^*$ variables that it adds to $\cS^*$ in line 6, all belong to $\{1,\ldots,k\}$, thus outputting $\cI^*$.

\vspace{0.1cm}

We begin by writing the distribution of the $i^{th}$ sample from  $\cN(0,\beta\vv^*{\vv^*}^T + I_p)$ explicitly as
\begin{equation}\label{eq:defOfDist}
    {\x}_i = \sbeta u_i\,{\vv^*} + \xii_i,
\end{equation}
where $\xii_i \in \R^p$ is a noise vector whose entries are all~\iid~$\cN(0,1)$, and $u_i\sim \cN(0,1)$. Furthermore, all the $u_i$'s and $\xii_i$'s are independent of each other.

\vspace{0.1cm}

By the greedy rule in line 3 of GreedySPCA, the $k-k^*$ variables in $\onep \setminus \cS^*$ that will be chosen are those with largest value of $f_{avg}(\cS^* \cup \{i \})$. We rewrite $f_{avg}(\cS^* \cup \{i \})$ as
\begin{align}\label{eq:ftilde}
    f_{avg}(\cS^* \cup \{i\})& = \tfrac{1}{k^*+1}\sum_{s,t \in \cS^* \cup \{i\}} \SigHat_{s,t} = \tfrac{k^*}{k^*+1} f_{avg}(\cS^*) + \tfrac{2}{k^*+1}\underbrace{ \left( \sum_{s \in \cS^*} \SigHat_{is}\right)}_{c_{i}(\cS^*)} +\tfrac{1}{k^*+1}\SigHat_{ii}:= \\& := \tfrac{k^*}{k^*+1} f_{avg}(\cS^*) + \nonumber \tfrac{1}{k^*+1}\left(2c_{i}(\cS^*)+\SigHat_{ii}\right).
\end{align}

The only part in Eq.~\eqref{eq:ftilde} that depends on $i$ is its total covariance with $\cS^*$ (which we denote by $c_i$), and its variance $\SigHat_{ii}$. In high level, the algorithm succeeds since $c_i$ is much bigger than $c_j$ for all pairs $i \in \cI^*,j \notin \cI^*$.
We now turn to make this argument formal.

\begin{lemma}\label{lem:ci_signal} Under the conditions of Theorem \ref{thm:main}, for a fixed $\cS^* \subseteq \cI^*$ of size $k^*$, w.p. at least $1-1/n$, every $i \in \cI^*$ satisfies $c_{i} \ge 0.4\beta k^*/k$.
\end{lemma}

\begin{lemma}\label{lem:ci_noise} Under the conditions of Theorem \ref{thm:main}, for a fixed $\cS^* \subseteq \cI^*$ of size $k^*$, w.p. at least $1-1/n$, every $j \notin \cI^*$ satisfies $c_{j} \le 0.3\beta k^*/k$.
\end{lemma}

\begin{lemma}\label{lem:DiagNoSig} Under the conditions of Theorem \ref{thm:main}, with probability at least $1-1/n$, for every $i \in \cI^*,j \notin \cI^*$,  $\SigHat_{jj}-\SigHat_{ii} \le 0.1\beta k^*/k.$
\end{lemma}

We  use Lemmas \ref{lem:ci_signal}, \ref{lem:ci_noise} and \ref{lem:DiagNoSig} to complete the proof of the theorem. Let $\Delta_{ij}=f_{avg}(\cS^{*} \cup \{i\})-f_{avg}(\cS^{*} \cup \{j\})$. To prove that GreedySPCA outputs $\cI^*$ we need to show that $\Delta_{ij}>0$ for every $i \in \cI^*,j \notin \cI^{*}$. Using Lemmas \ref{lem:ci_signal}--\ref{lem:DiagNoSig}, we have  $k \cdot \Delta_{ij} \ge 2(0.4\beta k^*/k-0.3\beta k^*/k)-0.1\beta k^*/k \ge 0.1\beta k^*/k >0$. In the last inequality we assumed $k^* \ge 1$.

The case $k^*=0$ corresponds to the regime $k=O(k_{comp}/\sqrt{\log p})$. In this regime the $k$ largest diagonal entries belong to $\cI^*$ \cite{Johnstone.Lu2009Consistency}. Indeed, when $k^*=0$ then $f_{avg}(\{i\})$ is simply $\SigHat_{ii}$, and GreedySPCA is no other than Diagonal Thresholding.

\vspace{0.1cm}

We  turn to prove Lemmas \ref{lem:ci_signal}--\ref{lem:DiagNoSig}. In the proof we use the following two auxiliary facts. The first is a large deviation result for a Chi-square random variable.
\begin{lemma}(\cite{LAURENT-MASSART90})\label{lem:XiSquare}
Let $X \sim \chi^2_{n}$. For all $x \geq 0$,
$$Pr[X \geq n+2\sqrt{nx} + x] \leq e^{-x},$$
$$Pr[X \leq n-2\sqrt{nx}] \leq e^{-x}.$$
\end{lemma}
The second fact records a well-known argument about the inner-product of two multivariate Gaussians.
\begin{lemma}\label{lem:CLT}
Let $\{x_i,y_i\}_{i=1}^n$  be standard~\iid~Gaussian random variables.
Then $\sum_{i=1}^n x_i y_i$ is distributed like the product
of two independent random variables $\|\x\|\cdot \tilde y$, where $\x=(x_1,\ldots,x_n)$,
$\|\x\|^2 \sim \chi^2_n$ and $\tilde y\sim\cN(0,1)$.
\end{lemma}
\begin{Proof}
For every fixed realization of $\x$,
we have $x_i y_i \sim \cN(0,x_i^2)$ and by the independence of the $y_i$'s,
$$
  \sum_{i=1}^n x_iy_i \sim  \cN(0,\|\x\|^2)=\|\x\|\cdot \cN(0,1):=\|\x\|\cdot \tilde y.$$
The lemma follows by observing that $\|\x\|^2 \sim \chi^2_n$.
\end{Proof}

\subsection{Proof of Lemma \ref{lem:ci_signal}}
We start by explicitly writing $c_i(\cS^*)$ from Eq.~\eqref{eq:ftilde} for a fixed set $\cS^*$ of size $k^*$. Let $\rr^{(i)}$ denote the $i^{th}$ row of the $p \times n$ design matrix $X$. For every candidate $i \notin \cS^*$,
$$c_i(\cS^*)=\sum_{j \in \cS^*} \SigHat_{ij}=\tfrac{1}{n}\sum_{j \in \cS^*} \rr^{(i)}\cdot (\rr^{(j)})^T= \tfrac{1}{n} \rr^{(i)}\left(\sum_{j \in \cS^*}  \rr^{(j)}\right)^T.$$
Following the distribution rule of $X$ given in Eq.~\eqref{eq:defOfDist}, all entries of the vector $\es=\sum_{j \in \cS^*}  \rr^{(j)}$ are i.i.d. with
$$\es_\ell \sim\sbeta u_\ell \sum_{j \in \cS^*}\vv^*_j+\sqrt{k^*}w_\ell,$$
where $u_\ell\sim \cN(0,1)$ is defined in Eq.~\eqref{eq:defOfDist} and $w_\ell  \sim \cN(0,1)$ independently of $u_\ell$ ($\sqrt{k^*}w_\ell$ is derived from $\sum_{j \in \cS^*} (\xii_\ell)_j$).
The product $\rr^{(i)}\es^T$ is distributed as
\begin{equation}\label{eq:distOfCov}
\rr^{(i)}\es^T \sim \frac{1}{n}\sum_{\ell=1}^n \left(\sbeta u_\ell \vv^*_i + y_\ell\right)\left(\sbeta u_\ell \left(\sum_{j \in \cS^*}\vv^*_j\right)+\sqrt{k^*}w_\ell\right)
\end{equation}
The variable $y_\ell=(\xii_\ell)_i \sim \cN(0,1)$. We rearrange the sum as four components, corresponding to the pure signal part, cross noise-signal and pure noise:
\begin{align}
&\sum_{\ell=1}^n \sbeta u_\ell \vv^*_i \cdot \sbeta u_\ell \sum_{j \in \cS^*}\vv^*_j = \frac{\beta k^* \vv^*_i}{\sqrt{k}}\sum_{\ell=1}^n u_\ell^2,\label{eq:onlysignal}
\\& \sum_{\ell=1}^n \sbeta u_\ell \vv^*_i \cdot \sqrt{k^*}w_\ell=  \sqrt{\beta k^*}\vv^*_i \sum_{\ell=1}^n u_\ell w_\ell,\label{eq:cross1}
\\& \sum_{\ell=1}^n y_\ell \sbeta u_\ell \sum_{j \in \cS^*}\vv^*_j =\frac{\sbeta k^*}{\sqrt{k}}\sum_{\ell=1}^n y_\ell u_\ell,\label{eq:cross2}
\\& \sum_{\ell=1}^n y_\ell \sqrt{k^*}  w_\ell=\sqrt{k^*}\sum_{\ell=1}^n y_\ell   w_\ell.\label{eq:noise}
\end{align}


We now bound each term separately. To lower bound Eq.~\eqref{eq:onlysignal}, we use the fact that  $\sum_{\ell=1}^n u_\ell^2$ in Eq.~\eqref{eq:onlysignal} is distributed $\chi^2_n$. The second inequality in Lemma \ref{lem:XiSquare} with $x=0.05n$ gives $Pr[\chi^2_n \le 0.8n]\le e^{-n/100}$. Therefore, w.p. at least $1-e^{-n/100}$,
\begin{equation}\label{eq:signalBound}
    \tfrac{1}{n}\eqref{eq:onlysignal} \ge 0.8\beta k^*/k
\end{equation}


Moving to Eq.~\eqref{eq:cross1}, according to Lemma \ref{lem:CLT}, the product term in Eq.~\eqref{eq:cross1} is distributed as $\sqrt{\chi^2_n}\cN(0,1)$. For $\tfrac{1}{n}\eqref{eq:cross1} > 0.1\beta k^*/k$ to hold,  the following has to happen,
$$\sqrt{\chi^2_n}|\cN(0,1)|>\frac{\beta n\sqrt{k^*}}{10k}=\frac{\beta n \sqrt{k^*}}{30k\sqrt{\log n}}\cdot \sqrt{9\log n}.$$
Using standard tail-bounds for Gaussians,
$Pr[|\cN(0,1)|>\sqrt{9\log n}] \le n^{-4}.$
Next we bound $Pr[\chi^2_n \ge \beta^2 n^2k^*/(900k^2\log n)].$ Substituting the value of $k^*$ from Eq.~\eqref{eq:kstar} we have that
$$\frac{\beta ^2 n^2k^*}{900k^2\log n} \ge \frac{\beta^2 n^2 }{ 900 k^2 \log n}\cdot \frac{Ck^2 \log n}{\beta^2 n}=Cn/900.$$
Choosing $C \ge 1800$ for example and using Lemma \ref{lem:CLT} gives $Pr[\chi^2_n \ge 2n]\le e^{-n/4}.$ To conclude, w.p. at least $1-n^{-4}-e^{-n/4}$ we get

\begin{equation}\label{eq:cross1bound}
   \tfrac{1}{n}|\eqref{eq:cross1} |\le \frac{0.1\beta k^*}{k}
\end{equation}


Moving to Eq.~\eqref{eq:cross2}, according to Lemma \ref{lem:CLT}, the sum-product term in Eq.~\eqref{eq:cross2} is distributed as $\sqrt{\chi^2_n}\cN(0,1)$.
Using  standard tail-bounds for Gaussians,
$Pr[|\cN(0,1)|\ge \sqrt{6\log n}]\le 2n^{-3},$
and according to Lemma \ref{lem:XiSquare}, $Pr[\chi^2_n \ge 2n]\le e^{-n/4}.$
Therefore w.p. at least $1-2n^{-3}-e^{-n/4}$ we get
\begin{equation}\label{eq:cross2bound}
    \tfrac{1}{n}\eqref{eq:cross2} \le \frac{1}{n}\frac{\sbeta k^*}{\sqrt{k}} \sqrt{6\log n}\sqrt{2n} \le \frac{0.1\beta k^*}{k}
\end{equation}
The last inequality is true when $k \le  \beta n / (1200\log n)$, which holds by our choice of $k$.

\vspace{0.1cm}

Moving to Eq.~\eqref{eq:noise}, we similarly have that w.p. at least $1-2n^{-3}-e^{-n/4}$
\begin{equation}\label{eq:noiusebound}
    \tfrac{1}{n}\eqref{eq:noise} \le \frac{1}{n}\sqrt{k^*}\sqrt{6\log n}\sqrt{2n} \le \frac{0.2k^*\beta}{k}.
\end{equation}
The last inequality in Eq.~\eqref{eq:noiusebound} holds whenever $k^* \ge 200k^2 \log n / (\beta^2 n)$, which is what Eq.~\eqref{eq:kstar} says.

Finally, the lower bound on $k^*$ in Eq.~\eqref{eq:kstar} makes sense as long as $k^* \le k$, which translates to requiring  $k \le \beta^2 n/(C \log n)$.

To conclude,  w.p. at least $1-3n^{-3}$, for a fixed $i \in \cI^*$,
$$c_i \ge \eqref{eq:signalBound} - \eqref{eq:cross1bound} - \eqref{eq:cross2bound}-\eqref{eq:noiusebound} \ge \frac{0.4k^*\beta}{k}.$$

The lemma now follows from taking the union bound over the $k-k^* \le p$ indices in $\cI^* \setminus \cS^*$, together with the fact that $p=O(n)$.

\subsection{Proof of Lemma  \ref{lem:ci_noise}}\label{sec:appdx:ci_noise}
For $i \notin \cI^*$ the terms in Eq.~\eqref{eq:onlysignal} and \eqref{eq:cross1} are 0 (because $\vv_i^*=0$). The proof of Lemma \ref{lem:ci_noise} is identical to the proof leading to the bound on Eq.~\eqref{eq:cross2} given in Eq.~\eqref{eq:cross2bound} and to the bound on Eq.~\eqref{eq:noise} given in Eq.~\eqref{eq:noiusebound}. A union bound is then taken over the at most $p$ variables in $\onep \setminus \cI^*$.

\subsection{Proof of Lemma  \ref{lem:DiagNoSig}}\label{sec:appdx:DiagNoSig}
For $j \notin \cI^*$, the distribution of $\SigHat_{jj} \sim \tfrac{\chi^2_n}{n}$ (From Eq.~\eqref{eq:distOfCov}).  Lemma \ref{lem:XiSquare} entails that for a fixed $j$, w.p. at least $1-n^{-3}$, $\SigHat_{jj}\le 1+\sqrt{\tfrac{9\log n}{n}}$.\\
For $i \in \cI^*$, $$\SigHat_{ii} \sim \frac{\beta}{k}\frac{\chi^2_n}{n}+2\frac{\sbeta}{\sqrt{k}}\frac{
\cN(0,1)\sqrt{\chi^2_n}}{n} +\frac{\chi^2_n}{n}$$
Using Lemma \ref{lem:XiSquare} and standard tails on the Gaussian, we obtain that w.p. at least $1-O(n^{-3})$, $\SigHat_{ii} \ge 1 + \tfrac{\beta}{k} - \sqrt{\tfrac{36\log n}{n}}.$
Using the union bound we get that w.p. at least $1-O(n^{-1})$, for every pair $i\in \cI^*,j \notin \cI^*$, $$\SigHat_{jj}-\SigHat_{ii} \le  \sqrt{\frac{100\log n}{n} }- \frac{\beta}{k}\le \sqrt{\frac{100\log n}{n} }  \le \frac{0.1\beta k^*}{k}.$$
The last inequality holds if $k^* \ge \sqrt{L}$, where $L$ is the lower bound on $k^*$ in Eq.~\eqref{eq:kstar}. However, $k^* \ge L$ implies $k^* \ge \sqrt{L}$ since $k^*$ is an integer.


\section{Proof of Theorem \ref{thm:spectral}}\label{sec:ProofSpectral}
We prove that w.p. tending to 1 as $(n,p,k) \to \infty$, for every set $\cI \subseteq \{1,\ldots,p\}$ of size $k$ that satisfies $| \cI \cap \cI^*|=\delta k$, for every $\delta \in [0,1]$:
\begin{equation}\label{eq:bound_lambda1}
\lambda_1(\SigHat_S) \in [1+\delta\beta -(1+2\sbeta)\Phi, 1+\delta\beta + (1+2\sbeta)\Phi +\frac{\beta}{k}],
\end{equation}
Where $\Phi=\sqrt{\frac{8  k \log n}{n}}$.

Fix a set $\cI \subseteq \{1,\ldots,p\}$ s.t. $|\cI \cap \cI^*|= \delta k$. The matrix $\SigHat_\cI$ can be written as $\SigHat_\cI = N + S$ where $N$ is composed of the noise part, Eq.~\eqref{eq:noise}, and $S$ is composed of the signal and noise-signal cross terms, Eq.~\eqref{eq:onlysignal}--\eqref{eq:cross2}.
$N$ is easily seen to be symmetric (in fact it follows a Wishart distribution), and therefore the matrix $S=\SigHat_\cI-N$, the difference of two symmetric matrices, is symmetric as well. Weyl's inequality, applicable for Hermitian matrices, implies that
\begin{equation}\label{eq:low_up_eig}
  \lambda_k(N)+\lambda_1(S) \le \lambda_1(\SigHat_\cI)\le \lambda_1(N)+\lambda_1(S)
\end{equation}
\subsection{Bounding $\lambda_1(N)$ and $\lambda_k(N)$}
The matrix $N \in \R^{k \times k}$ follows a Wishart distribution, and  by~\cite[Theorem II.13]{Szarek:survey},
\begin{equation*}
\pr[\lambda_1(N) \ge (1+\sqrt{k/n}+t)^2 \,\vee\, \lambda_k(N) \le (1-\sqrt{k/n}-t)^2] \le e^{-nt^2/2}.
\end{equation*}
Plugging in $t=\sqrt{6k \log n/n}$ we obtain that w.p. at least $1-n^{-3k}$,
\begin{equation} \label{eq:wishart_eig_bound1}
\lambda_1(N) \le \left(1+\sqrt{\frac{k}{n}}+\sqrt{\frac{6k\log n}{n}}\right)^2 \le 1+\sqrt{\frac{8k \log n}{n}} = 1+\Phi.
\end{equation}
and similarly,
\begin{equation} \label{eq:wishart_eig_bound2}
\lambda_k(N) \ge  1-\Phi.
\end{equation}
Taking the union bound over all $\binom{p}{k} \le p^k$ possible sub-matrices $N$, the bounds hold w.p. at least $1-n^{-3k}p^k \ge 1-n^{-1}$.

\subsection{Upper Bounding $\lambda_1(S)$}
Recall the parameters $\uu=(u_1,\ldots,u_n)$ and $\xii_i$ from definition of the single-spike distribution given in Eq.~\eqref{eq:defOfDist}.
We start with a certain property of $\SigHat$ that we require during the proof of the upper and lower bound on $\lambda_1(S)$. We say that $\SigHat$ is {\it{typical}} if $\| \uu \|^2 \le n+6\sqrt{n\log n} $ and if $(\xii_1)_i \le 2\sqrt{\log n}$  for every $i=1,\ldots,p$. Lemmas \ref{lem:XiSquare} and \ref{lem:CLT} guarantee that $\SigHat$ is typical w.p. at least $1-n^{-1}$. In what follows we condition on this fact.

To upper bound the largest eigenvalue of $S$ we use Gershgorin's circle theorem, which says that every eigenvalue $\lambda$ of an $n \times n$ matrix $A$ satisfies at least one of the $n$ inequalities for $i=1,\ldots,n$,
\begin{equation}\label{eq:gershgo}
|\lambda-A_{ii}| \le \sum_{j \ne i} |A_{ij}|.
\end{equation}
Each inequality defines a Gershgorin's disc, and every $\lambda$ belongs to at least one disc. We next show that all discs are almost identical, and evaluate their center and radius.

Decompose each entry $S_{ij}$ according to the three sums Eq.~\eqref{eq:onlysignal}-\eqref{eq:cross2} (plugging $k^*=1$). To bound the sums in Eq.\eqref{eq:cross1} and \eqref{eq:cross2} we note that both involve the term $u_\ell$, which does not depend on $i$ or $j$. Therefore we may rotate the distribution to point in the direction of $\uu$. According to Lemma \ref{lem:CLT},  the sum-product \eqref{eq:cross1}  is then distributed $\|\uu\|\cdot (\xii_1)_i$ and Eq.\eqref{eq:cross2}  is  distributed $\|\uu\|\cdot (\xii_1)_j$.
Using the fact that $\SigHat$ is typical we obtain the following bounds:
$$\tfrac{1}{n}|\eqref{eq:cross1}| \sim \tfrac{1}{n}\sbeta\vv_i^*\|\uu\| |(\xii_1)_j| \le \sqrt{\frac{8\beta  \log n}{nk}}.$$
Similarly $$\tfrac{1}{n}|\eqref{eq:cross2}| \le \sqrt{\frac{8\beta  \log n}{nk}}, \qquad \tfrac{1}{n}|\eqref{eq:onlysignal}|=\left(1 \pm \sqrt{\frac{36 \log n}{n}}\right)\frac{\beta}{k} .$$

Putting everything together, and letting $\delta_i=1$ if $i \in \cI^*$ and 0 otherwise, if $\SigHat$ is typical then for every $i,j\in \cI$,
\begin{equation}\label{eq:Sij}
    S_{ij} = \delta_i\delta_j\frac{\beta}{k}+\Delta_{ij}, \qquad |\Delta_{ij}| \le \Delta := \sqrt{\frac{36 \beta \log n}{nk}}.
\end{equation}
To bound the radius of the $i^{th}$ disc, $\sum_j |S_{ij}|$, we need to account for $|\cI \cap \cI^*|=\delta k$ indices $j \in \cI^*$ and $(1-\delta)k$ indices $j \notin \cI^*$. Plugging \eqref{eq:Sij} in \eqref{eq:gershgo}, we obtain that
$$|\lambda - S_{ii}| \le \delta k  \left(\frac{\beta}{k}+\Delta\right) + (1-\delta) k \cdot \Delta=\delta\beta+\Delta k.$$
Rearranging we get
\begin{equation}\label{eq:upper_bound_S}
\lambda_1(S) \le S_{ii} + \delta \beta +\Delta k \le \frac{\beta}{k}+\Delta+\delta \beta +\Delta k \le \delta \beta + \left(\frac{\beta}{k}+2\Delta k \right).
\end{equation}

\subsection{Lower Bounding $\lambda_1(S)$}
To lower bound the largest eigenvalue of $S$ we use the Rayleigh quotient definition, namely $\lambda_1(S)$ is the argmax of   $\x^T S \x$ over all unit vectors $\x \in \R^k$. In particular, for  $\x_0=(\delta k)^{-0.5}\bf{1}_{\cI \cap \cI^*}$ (${\bf{1}}_Q$ is the characteristic vector of a set $Q$), the value of  $\x_0^T S \x_0$ is  a lower bound on $\lambda_1(S)$. The latter is simply the average row sum in the $\delta k \times \delta k$ submatrix  $S_{\cI \cap \cI^*}$. If $\SigHat$ is typical, then according to Eq.~\eqref{eq:Sij},
\begin{equation}\label{eq:lower_bound_S}
\lambda_1(S) \ge  \delta k \cdot \left(\frac{\beta}{k}-\Delta\right) \ge \delta \beta - \Delta k.
\end{equation}

To conclude the proof of the theorem, note that $\Delta k \le 3 \sbeta\Phi$. Putting Equations \eqref{eq:wishart_eig_bound1},\eqref{eq:wishart_eig_bound2},\eqref{eq:upper_bound_S},\eqref{eq:lower_bound_S} together, we get that  w.p. at least $1-2n^{-1}$,
$$\lambda_1(\SigHat_S) \in [1+\delta\beta -(1+3\sbeta)\Phi, 1+\delta\beta  + (1+3\sbeta)\Phi+\frac{\beta}{k}].$$

\subsection{Proof of Corollary \ref{cor:xi}}\label{sec:proofOfCor}
Take $\cI,\cJ \subseteq \onep$ that satisfy $|\cI \cap \cI^*|=\delta_1 > \delta_2 = |\cJ \cap \cI^*|$. According to Theorem \ref{thm:spectral}, for $\lambda_1(\SigHat_{\cI})>\lambda_1(\SigHat_{\cJ})$ to hold, it suffices to require
$$ 1+\delta_2\beta+\frac{\beta}{k} + \Gamma < 1+\delta_1\beta -\Gamma.$$
Rearranging we get,
$$\frac{2\Gamma}{\beta} + \frac{1}{k} < \delta_1 - \delta_2 := \xi.$$
The corollary follows immediately from the definition of $\Gamma$ and $k_{info}$.


\section{Simulations}\label{sec:simulation}
In this section, we evaluate the performance of SSPCA both in the strong and weak SNR regimes.
The definition of  $k_{comp}$ in Section \ref{sec:spikedmodel} (the boundary between the two  regimes) is asymptotic and cannot be used directly in simulation. Instead, we let the empirical success rate of Covariance Thresholding (abbreviated CT) define regime boundaries. The following points summarize the way we ran simulations:
\begin{itemize}
    \item  We sample from the uniform \emph{unbiased} sparse PCA distribution (UNBSPCA), where $\vv^*=\left(\pm \frac{1}{\sqrt{k}},\ldots,\pm\frac{1}{\sqrt{k}},0,\ldots,0\right)$. The signs of non-zero entries are randomly chosen. Accordingly, we change the choice of $f_1$ in GreedySPCA from $f_{avg}$ to $f_{\ell_1}$, which measures the $\ell_1$ norm of $\SigHat_{S^* \cup \{i\}}$, rather than the sum. Furthermore, in the weak SNR regime, it makes sense to ignore the diagonal of $\SigHat$. Therefore we only use $c_i$ from Eq.~\eqref{eq:ftilde} to choose $i$, with $\ell_1$ norm.

    \item We keep $p=n$ and fix $\beta=0.5$. The choice of $0.5$ is somewhat arbitrary and any value below $\sqrt{p/n}=1$ is suitable. When $\beta$ exceeds $\sqrt{p/n}$ the problem is computationally easy for all  $k$ up to the information limit~\cite{krauthgamer2015}[Thm 1.1].
    \item The \emph{success rate} of an algorithm on a given input is defined to be $\tfrac{1}{k}|\cI \cap \cI^*|$, where $\cI\subseteq\{1,\ldots,p\}$ is the algorithm's guess of $\vv^*$'s support.
    \item The algorithm DT has no tunable parameters -- it simply returns the indices of the $k$ largest diagonal entries. The performance of CT, on the other hand, depends crucially on the chosen threshold. When running CT we loop over 50 thresholds, the empirical percentiles of the off-diagonal entries of the input covariance matrix. We choose the best result as CT's output.
    \item Formally, the output of CT is a vector (a guess for $\vv^*$). We convert the vector to a set $\cI$ by taking the indices of the $k$ largest entries in absolute value.
    \item We  compared SSPCA against the well-known sparse PCA algorithm described in \cite{Zou04sparseprincipal}. This algorithm casts PCA as a regression problem and uses both ridge and lasso penalties. LARS~\cite{efron2004} is then used to obtain the optimal solution. We ran Python's implementation of this algorithm using a grid search for the ridge and lasso penalties in the rectangle $[0,2]\times [0,2]$, discretized to 100 equidistant points. The algorithm is in module \texttt{sklearn.decomposition.SparsePCA}~\cite{scikit-learn}.
\end{itemize}

Figure \ref{fig:comparison} compares the performance of all aforementioned algorithms in a certain weak SNR configuration, $n=p=1000,k=8,\beta=0.5$. Among the polynomial-time algorithms (we include in this category SSPCA with $k^*=1,2$), SSPCA with $k^*=2$ performs best. When running for super polynomial-time, SSPCA with $k^*=3$ is superior to  the naive exhaustive search, when both are given the same time budget $T$ (3 parallel-hours on 90 cores).
\begin{figure}[!t]
\centering
\adjustimage{width=19cm,center}{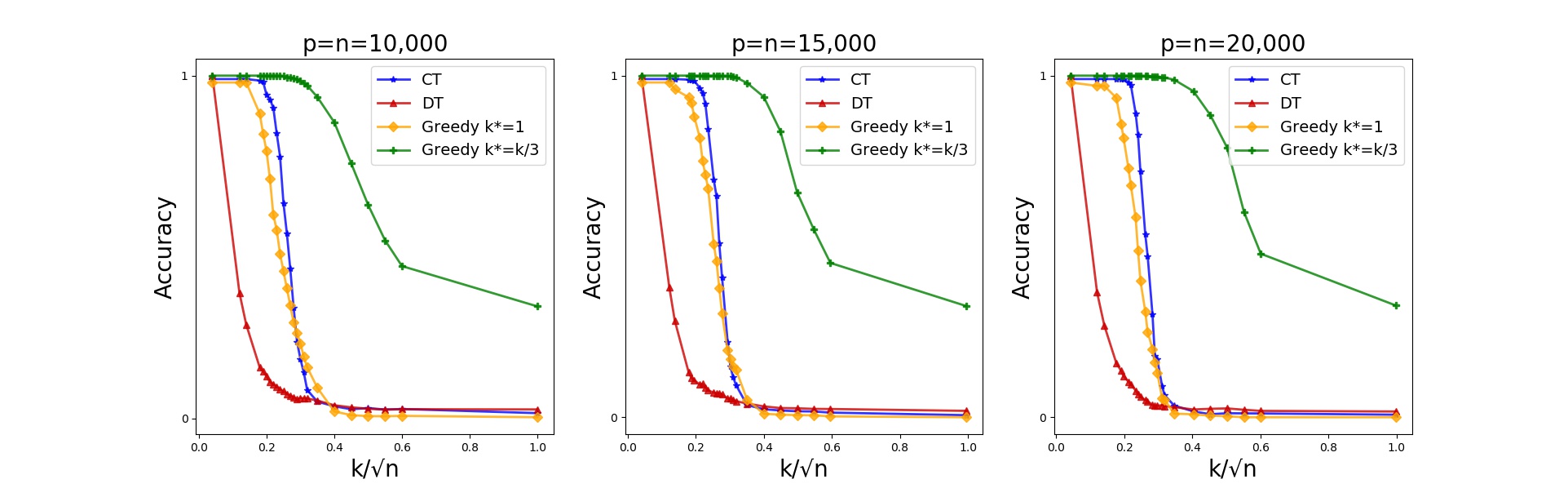}
\caption{The success rate of DT, CT and GreedySPCA as a function of $k$. Every point is an average of 25 executions, with $n=p$ samples. GreedySPCA was initialized with $k^*=1$ or $k^*=k/3$ random entries from $\cI^*$.}
\label{fig:boundaries}
\end{figure}

Our next experiment demonstrates the existence of golden seeds in the weak SNR regime. The empirical boundary of the strong/weak SNR regime is charted by the success rate curve of CT. Figure \ref{fig:boundaries} shows the performance ($y$-axis) of CT as $k$ increases. Three configurations are plotted $n=p=10,000,15,000,20,000$. The $x$-axis is scaled by $\sqrt{n}$ to defuse the dependence on $n$. Indeed all three CT lines overlap as expected (due to scaling), and the phase transition to the hard regime occurs when $k$ is in the window $[0.2\sqrt{n},0.3\sqrt{n}]$. The plot also includes the performance of DT, lagging behind, and GreedySPCA initialized  once with a seed of size $k^*=1$ and second with $k^*=k/3$. In both cases the seed is a random subset of $\cI^*$.

\begin{figure}[!t]
\centering
\adjustimage{width=6cm,center}{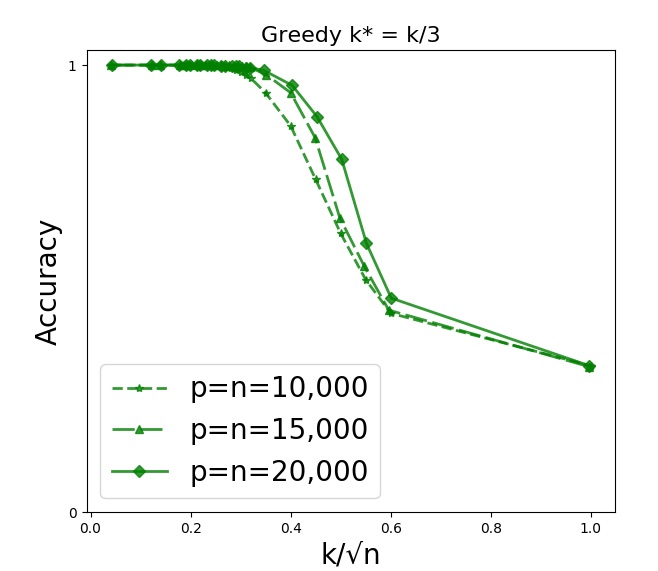}
\caption{The success rate of  GreedySPCA when initialized with  $k^*=k/3$ random entries from $\cI^*$. Lines are non-overlapping due to under-scaling on the $x$-axis.}
\label{fig:weakSNR}
\end{figure}

As evident from Figure \ref{fig:boundaries}, the performance of GreedySPCA with seeds of size $k^*=1$ is similar to CT. This is somewhat surprising when comparing to the asymptotic lower bound given by Eq.~\eqref{eq:kstar}, which would be of order $\log n$ at the computational threshold $k_{comp}$. The right-most (green) line in Figure~\ref{fig:boundaries} extends with an accuracy of roughly $100\%$ into the weak SNR regime, thus showing the existence of golden seeds in that regime. Another cue that we are indeed in the weak SNR regime is the fact that the lines in Figure \ref{fig:weakSNR} corresponding to  different input sizes do not overlap, suggesting that the scaling of the $x$-axis is too small.

\section{Discussion}
In this paper, we presented a family of anytime algorithms for the $k$-sparse PCA problem, which follow the same simple white-box greedy template that we called  GreedySPCA and SSPCA.

GreedySPCA performs a bulk greedy choice, and instead we could have grown the solution iteratively, adding in iteration $r=1,\ldots,k-k^*$ the variable $i_r$ which maximizes $f_1(\cS^* \cup \{i_1,\ldots,i_{r-1}\})$. The iterative variant is exactly the well-known greedy algorithm of Nemhauser, Wolsey and Fisher, which was proposed for sub-modular function optimization~\cite{Nemhauser:1978}. The only difference is that Nemhauser et al. start with an empty seed.

Nemhauser et al. proved that if $f_1$ is sub-modular and monotone, then the iterative greedy algorithm finds a solution which is a $(1-\tfrac{1}{e})$-approximation of the optimum. Hence it is tempting to run the iterative greedy with $f_1=f_{\lambda_1}$ and no seed. However, even if  $f_{\lambda_1}$ is monotone and sub-modular, the $(1-\tfrac{1}{e})$-approximation ratio is useless in many cases, as the guaranteed value is lower than a random solution (see Theorem \ref{thm:spectral}). Our proof shows that GreedySPCA recovers $\cI^*$ exactly when called with the right seed without the sub-modularity assumption on $f_1$.

Simulations that we ran in the spiked covariance model with the iterative version (with seed) resulted in very similar performance compared to the bulk version. Similarly, using $f_{\lambda_1}$ instead of $f_{\ell_1}$ in GreedySPCA did not incur  improvement.

\section*{Acknowledgment} We thank Jonathan Rosenblatt for allowing us to use his cluster, and for his patient technical support.

\bibliographystyle{alphaurlinit}
\bibliography{bibfile}

\newcommand{\etalchar}[1]{$^{#1}$}
\begin{thebibliography}{dEGJL04}

\bibitem[And84]{Anderson84}
T.~Anderson.
\newblock {\em An introduction to multivariate statistical analysis}.
\newblock Wiley series in probability and mathematical statistics. Wiley, 2nd
  edition, 1984.

\bibitem[AW09]{AminiWain09}
A.~Amini and M.~Wainwright.
\newblock High dimensional analysis of semidefinite relaxations for sparse
  principal component analysis.
\newblock {\em Annals of Statistics}, 37(5B):2877--2921, 2009.
\newblock \href {http://dx.doi.org/10.1214/08-AOS664}
  {\path{doi:10.1214/08-AOS664}}.

\bibitem[BB19]{Bresler}
M.~Brennan and G.~Bresler.
\newblock Optimal average-case reductions to sparse pca: From weak assumptions
  to strong hardness.
\newblock In {\em COLT}, 02 2019.

\bibitem[BBH18]{pmlr-v75-brennan18a}
M.~Brennan, G.~Bresler, and W.~Huleihel.
\newblock Reducibility and computational lower bounds for problems with planted
  sparse structure.
\newblock In {\em Proceedings of the 31st Conference On Learning Theory},
  volume~75 of {\em Proceedings of Machine Learning Research}, pages 48--166.
  PMLR, 06--09 Jul 2018.

\bibitem[BL08]{BickelLevina06}
J.~Bickel and E.~Levina.
\newblock Regularized estimation of large covariance matrices.
\newblock {\em Annals of Statistics}, 36:199--227, 2008.
\newblock \href {http://dx.doi.org/10.1214/009053607000000758}
  {\path{doi:10.1214/009053607000000758}}.

\bibitem[BR13a]{RigolletCliqueCOLT}
Q.~Berthet and P.~Rigollet.
\newblock Complexity theoretic lower bounds for sparse principal component
  detection.
\newblock In {\em COLT}, pages 1046--1066, 2013.
\newblock \href {http://arxiv.org/abs/1304.0828} {\path{arXiv:1304.0828}}.

\bibitem[BR13b]{RigolletClique}
Q.~Berthet and P.~Rigollet.
\newblock Optimal detection of sparse principal components in high dimension.
\newblock {\em Annals of Statistics}, 41(4):1780--1815, 08 2013.
\newblock \href {http://dx.doi.org/10.1214/13-AOS1127}
  {\path{doi:10.1214/13-AOS1127}}.

\bibitem[CMW13]{CaiMaWu13}
T.~Cai, Z.~Ma, and Y.~Wu.
\newblock {Sparse PCA: Optimal rates and adaptive estimation}.
\newblock {\em The Annals of Statistics}, 41(6):3074--3110, 2013.
\newblock \href {http://dx.doi.org/10.1214/13-AOS1178}
  {\path{doi:10.1214/13-AOS1178}}.

\bibitem[CMW15]{CMW15}
T.~Cai, Z.~Ma, and Y.~Wu.
\newblock Optimal estimation and rank detection for sparse spiked covariance
  matrices.
\newblock {\em Probability Theory and Related Fields}, 161(3-4):781--815, 4
  2015.
\newblock \href {http://dx.doi.org/10.1007/s00440-014-0562-z}
  {\path{doi:10.1007/s00440-014-0562-z}}.

\bibitem[dEGJL04]{AspremontEtAlSDP07}
A.~d'Aspremont, L.~El-Ghaoui, M.~Jordan, and G.~Lanckriet.
\newblock A direct formulation for sparse {PCA} using semidefinite programming.
\newblock {\em SIAM Review}, 49(3):434--448, 2004.
\newblock \href {http://dx.doi.org/10.1137/050645506}
  {\path{doi:10.1137/050645506}}.

\bibitem[DKWB19]{ding2019subexponentialtime}
Y.~Ding, D.~Kunisky, A.~S. Wein, and A.~S. Bandeira.
\newblock Subexponential-time algorithms for sparse pca, 2019.
\newblock \href {http://arxiv.org/abs/1907.11635} {\path{arXiv:1907.11635}}.


\bibitem[DM16]{Deshpande:2016}
Y.~Deshpande and A.~Montanari.
\newblock Sparse pca via covariance thresholding.
\newblock {\em J. Mach. Learn. Res.}, 17(1):4913--4953, January 2016.
\newblock Available from:
  \url{http://dl.acm.org/citation.cfm?id=2946645.3007094}.

\bibitem[Don00]{Donoho00}
D.~Donoho.
\newblock High-dimensional data analysis: The curses and blessings of
  dimensionality.
\newblock In {\em AMS CONFERENCE ON MATH CHALLENGES OF THE 21ST CENTURY}, 2000.

\bibitem[DS01]{Szarek:survey}
K.~Davidson and S.~Szarek.
\newblock Local operator theory, random matrices and {B}anach spaces.
\newblock In Lindenstrauss, editor, {\em Handbook on the Geometry of Banach
  spaces}, volume~1, pages 317--366. Elsevier Science, 2001.

\bibitem[EHJT04]{efron2004}
B.~Efron, T.~Hastie, I.~Johnstone, and R.~Tibshirani.
\newblock Least angle regression.
\newblock {\em Annals of Statistics}, 32(2):407--499, 4 2004.
\newblock Available from: \url{https://doi.org/10.1214/009053604000000067},
  \href {http://dx.doi.org/10.1214/009053604000000067}
  {\path{doi:10.1214/009053604000000067}}.

\bibitem[JL09]{Johnstone.Lu2009Consistency}
I.~M. Johnstone and A.~Lu.
\newblock {On Consistency and Sparsity for Principal Components Analysis in
  High Dimensions}.
\newblock {\em Journal of the American Statistical Association},
  104(486):682--693, 2009.
\newblock \href {http://dx.doi.org/10.1198/jasa.2009.0121}
  {\path{doi:10.1198/jasa.2009.0121}}.

\bibitem[Joh01]{Johnstone01}
I.~M. Johnstone.
\newblock On the distribution of the largest eigenvalue in principal components
  analysis.
\newblock {\em Annals of Statistics}, 29:295--327, 2001.
\newblock \href {http://dx.doi.org/10.1214/aos/1009210544}
  {\path{doi:10.1214/aos/1009210544}}.

\bibitem[Jol02]{JoliffePCA}
I.~Jolliffe.
\newblock {\em Principal Component Analysis}.
\newblock Springer series in statistics. Springer, 2nd edition, 2002.

\bibitem[KNV15]{krauthgamer2015}
R.~Krauthgamer, B.~Nadler, and D.~Vilenchik.
\newblock Do semidefinite relaxations solve sparse pca up to the information
  limit?
\newblock {\em Ann. Statist.}, 43(3):1300--1322, 06 2015.
\newblock Available from: \url{https://doi.org/10.1214/15-AOS1310}, \href
  {http://dx.doi.org/10.1214/15-AOS1310} {\path{doi:10.1214/15-AOS1310}}.

\bibitem[LM00]{LAURENT-MASSART90}
B.~Laurent and P.~Massart.
\newblock Adaptive estimation of a quadratic functional by model selection.
\newblock {\em Annals of Statistics}, 28(5):1302--�1338, 2000.
\newblock \href {http://dx.doi.org/10.1214/aos/1015957395}
  {\path{doi:10.1214/aos/1015957395}}.

\bibitem[Mui82]{Muirhead1982}
R.~J. Muirhead.
\newblock {\em Aspects of Multivariate Statistical Theory}.
\newblock Wiley, New York, 1982.

\bibitem[Nad08]{Nadler08}
B.~Nadler.
\newblock Finite sample approximation results for principal component analysis:
  a matrix perturbation approach.
\newblock {\em Annals of Statistics}, 36:2791�--2817, 2008.
\newblock \href {http://dx.doi.org/10.1214/08-AOS618}
  {\path{doi:10.1214/08-AOS618}}.

\bibitem[Nat95]{Natarajan1995}
B.~K. Natarajan.
\newblock Sparse approximate solutions to linear systems.
\newblock {\em SIAM J. Comput.}, 24(2):227--234, 1995.

\bibitem[NWF78]{Nemhauser:1978}
G.~L. Nemhauser, L.~A. Wolsey, and M.~L. Fisher.
\newblock An analysis of approximations for maximizing submodular set
  functions--i.
\newblock {\em Math. Program.}, 14(1):265--294, December 1978.
\newblock Available from: \url{https://doi.org/10.1007/BF01588971}, \href
  {http://dx.doi.org/10.1007/BF01588971} {\path{doi:10.1007/BF01588971}}.

\bibitem[Pau07]{Paul07}
D.~Paul.
\newblock Asymptotics of sample eigenstructure for a large dimensional spiked
  covariance model.
\newblock {\em Statistica Sinica}, 17:1617�--1642, 2007.

\bibitem[PVG{\etalchar{+}}11]{scikit-learn}
F.~Pedregosa, G.~Varoquaux, A.~Gramfort, V.~Michel, B.~Thirion, O.~Grisel,
  M.~Blondel, P.~Prettenhofer, R.~Weiss, V.~Dubourg, J.~Vanderplas, A.~Passos,
  D.~Cournapeau, M.~Brucher, M.~Perrot, and E.~Duchesnay.
\newblock Scikit-learn: Machine learning in {P}ython.
\newblock {\em Journal of Machine Learning Research}, 12:2825--2830, 2011.

\bibitem[SSM13]{SHEN2013317}
D.~Shen, H.~Shen, and J.~Marron.
\newblock Consistency of sparse pca in high dimension, low sample size
  contexts.
\newblock {\em Journal of Multivariate Analysis}, 115:317 -- 333, 2013.
\newblock Available from:
  \url{http://www.sciencedirect.com/science/article/pii/S0047259X12002308},
  \href {http://dx.doi.org/https://doi.org/10.1016/j.jmva.2012.10.007}
  {\path{doi:https://doi.org/10.1016/j.jmva.2012.10.007}}.

\bibitem[VL12]{minmax-pmlr-v22-vu12}
V.~Vu and J.~Lei.
\newblock Minimax rates of estimation for sparse pca in high dimensions.
\newblock In {\em Proceedings of the Fifteenth International Conference on
  Artificial Intelligence and Statistics}, volume~22 of {\em Proceedings of
  Machine Learning Research}, pages 1278--1286, La Palma, Canary Islands,
  21--23 Apr 2012. PMLR.

\bibitem[WBS14]{Samworth14}
T.~Wang, Q.~Berthet, and R.~Samworth.
\newblock Statistical and computational trade-offs in estimation of sparse
  principal components.
\newblock \url{http://arxiv.org/abs/1408.5369}, August 2014.

\bibitem[WBS16]{wang2016}
T.~Wang, Q.~Berthet, and R.~Samworth.
\newblock Statistical and computational trade-offs in estimation of sparse
  principal components.
\newblock {\em Ann. Statist.}, 44(5):1896--1930, 10 2016.
\newblock Available from: \url{https://doi.org/10.1214/15-AOS1369}, \href
  {http://dx.doi.org/10.1214/15-AOS1369} {\path{doi:10.1214/15-AOS1369}}.

\bibitem[ZHT06]{Zou04sparseprincipal}
H.~Zou, T.~Hastie, and R.~Tibshirani.
\newblock Sparse principal component analysis.
\newblock {\em Journal of Computational and Graphical Statistics},
  15(2):265--286, 2006.
\newblock Available from: \url{https://doi.org/10.1198/106186006X113430}, \href
  {http://arxiv.org/abs/https://doi.org/10.1198/106186006X113430}
  {\path{arXiv:https://doi.org/10.1198/106186006X113430}}, \href
  {http://dx.doi.org/10.1198/106186006X113430}
  {\path{doi:10.1198/106186006X113430}}.

\bibitem[Zil96]{Zilberstein_1996}
S.~Zilberstein.
\newblock Using anytime algorithms in intelligent systems.
\newblock {\em AI Magazine}, 17(3):73, Mar. 1996.
\newblock Available from:
  \url{https://www.aaai.org/ojs/index.php/aimagazine/article/view/1232}, \href
  {http://dx.doi.org/10.1609/aimag.v17i3.1232}
  {\path{doi:10.1609/aimag.v17i3.1232}}.



\end{thebibliography}
\end{document}